\documentclass{amsart}

\usepackage{color,graphicx,amssymb,latexsym,amsfonts,txfonts,amsmath,amsthm}
\usepackage{amscd}
\usepackage{pdfsync}

\usepackage{hyperref}
\hypersetup{
    colorlinks=true,       % false: boxed links; true: colored links
    linkcolor=blue,          % color of internal links
    citecolor=blue,        % color of links to bibliography
    filecolor=blue,      % color of file links
    urlcolor=blue           % color of external links
}

\newtheorem{theo}{Theorem}                     
                     
\newtheorem{coro}{Corollary}

\theoremstyle{remark}                  
\newtheorem{rema}{\bf Remark}
\newtheorem{example}{\bf Example}

\begin{document}
\title{Homology group automorphisms of Riemann surfaces}
\author{Rub\'en A. Hidalgo}
\thanks{Partially supported by project Fondecyt 1190001}
\address{Departamento de Matem\'atica y Estad\'{\i}stica, Universidad de La Frontera. Temuco, Chile}
\email{ruben.hidalgo@ufrontera.cl}

\begin{abstract}
If $\Gamma$ is a finitely generated Fuchsian group such that its derived subgroup $\Gamma'$ is co-compact and torsion free, then $S={\mathbb H}^{2}/\Gamma'$ is a closed Riemann surface of genus $g \geq 2$ admitting the abelian group $A=\Gamma/\Gamma'$ as a group of conformal automorphisms. We say that $A$ is a homology group of $S$. A natural question is if $S$ admits unique homology groups or not, in other words, is there are different Fuchsian groups $\Gamma_{1}$ and $\Gamma_{2}$ with $\Gamma_{1}'=\Gamma'_{2}$? It is known that if $\Gamma_{1}$ and $\Gamma_{2}$ are both of the same signature $(0;k,\ldots,k)$, for some $k \geq 2$, then the equality $\Gamma_{1}'=\Gamma_{2}'$ ensures that $\Gamma_{1}=\Gamma_{2}$. Generalizing this, we observe that if $\Gamma_{j}$ has signature $(0;k_{j},\ldots,k_{j})$ and $\Gamma_{1}'=\Gamma'_{2}$, then $\Gamma_{1}=\Gamma_{2}$. We also provide examples 
of surfaces $S$ with different homology groups. A description of the normalizer in ${\rm Aut}(S)$
of each homology group $A$ is also obtained.
\end{abstract}

\maketitle

%%%%%%%%%%%%%%%%%%%%
\section{Introduction}
Let $S$ be a closed Riemann surface of genus $g \geq 2$ and let ${\rm Aut}(S)$ be its group of conformal automorphisms. In 1890, Schwarz \cite{Schwarz} proved that ${\rm Aut}(S)$ is finite and later, in 1893, Hurwitz \cite{Hurwitz} obtained the upper bound  $|{\rm Aut}(S)| \leq 84(\gamma-1)$. Since then, the study of groups of conformal automorphisms of closed Riemann surfaces has been of interest in the community of Riemann surfaces and related areas.
In 1987, Nakayima \cite{Nakayima} proved that if $A<{\rm Aut}(S)$ is an abelian group, then $|A| \leq 4(g+1)$ (if $A$ is a cyclic group, then $|A| \leq 4g+2$ \cite{Wiman}).

An abelian group $A<{\rm Aut}(S)$ is called a {\it homology group} of $S$ if there is not a closed Riemann surface $R$ of genus $\gamma>g$ admitting an abelian group $B$ of conformal automorphisms such that $R/B$ and $S/A$ are isomorphic as Riemann orbifolds. In this case, $S$ is a {\it homology Riemann surface} and $(S,A)$ is a {\it homology Riemann pair}.

In terms of Fuchsian groups, the above can be described as follows (see \cite{Hidalgo:Homology}). Let ${\mathbb H}^{2}$ be the hyperbolic plane, let $\Gamma<{\rm PSL}_{2}({\mathbb R})={\rm Aut}({\mathbb H}^{2})$ be a Fuchsian group such that $S/A={\mathbb H}^{2}/\Gamma$ (as Riemann orbifolds) and let $\Gamma'$ be it derived subgroup. Then $A$ is a homology group of $S$ if and only if: (i) $\Gamma'$ is torsion free, (ii) $S={\mathbb H}^{2}/\Gamma'$ and (iii) $A=\Gamma/\Gamma'$.

Note that, as the index of $\Gamma'$ in $\Gamma$ is finite (since $A$ is finite), conditions (ii) and (iii) necessarily assert that $S/A$ has genus zero, so  it  has signature of the form $(0;k_{1},\ldots,k_{n+1})$, for some $n \geq 2$ and $k_{j} \geq 2$. Condition (i)  is equivalent to satisfies  
Maclachlan's condition \cite{Maclachlan}
\begin{equation}\label{E1}
{\rm lcm}(k_{1},...,k_{j-1},k_{j+1},...,k_{n+1})={\rm lcm}(k_{1},...,k_{n+1}), \quad \forall j=1,...,n+1,
\end{equation}
where ${\rm lcm}$ denotes the ``least common multiple".

Both, (i) a description of those hyperelliptic homology Riemann surfaces and (ii) an algebraic representation of those homology Riemann pairs $(S,A)$ where $S/A$ has a triangular signature, were given in \cite{Hidalgo:Homology}. In the same paper it was noticed that a homology group  
cannot be isomorphic to ${\mathbb Z}_{p}$, for $p$ prime, nor to ${\mathbb Z}_{2}^{2}$. But it 
can be isomorphic to a cyclic group of order different from a prime. For instance, if $\Gamma=\langle x_{1},x_{2},x_{3}: x_{1}^{2}=x_{2}^{5}=x_{3}^{10}=x_{1}x_{2}x_{3}=1\rangle$, then  $A=\Gamma/\Gamma' \cong {\mathbb Z}_{10}$ is a homology group for the homology Riemann surface $S={\mathbb H}^{2}/\Gamma'$.

\medskip

Let $(S={\mathbb H}^{2}/\Gamma',A=\Gamma/\Gamma')$ be a homology Riemann pair such that 
all the cone points of $S/A$ have the same order $k \geq 2$, i.e., $\Gamma$ has signature $(0;k,\stackrel{n+1}{\ldots},k)$ for some $k,n \geq 2$, 
then $A \cong {\mathbb Z}_{k}^{n}$. In this particular situation, we say that $A$ is a {\it generalized Fermat group of type $(k,n)$} and that $S$ is a {\it generalized Fermat curve of type $(k,n)$}.
By the Riemann-Hurwitz formula, the genus of $S$ is $g=1+k^{n-1}((n-1)(k-1)-2)/2$, so (i)  $(k-1)(n-1)>2$ and (ii) $n$ is uniquely determined by $g$ and $k$. An algebraic curve description for $S$ was obtained in \cite{GHL} (see Section \ref{Sec:GFC}). 

In \cite{HLKP} it was proved that a homology Riemann surface $S$ admits at most one generalized Fermat group of a fixed type $(k,n)$, equivalently, if $\Gamma_{1}$ and $\Gamma_{2}$ are both Fuchsian groups with the same signature $(0;k,\stackrel{n+1}{\ldots},k)$ such that $\Gamma'_{1}=\Gamma'_{2}$, then $\Gamma_{1}=\Gamma_{2}$.

This does not rule out the possibility for $S$ to have generalized Fermat groups of different types.
We start by observing that this is not the situation.

\begin{theo}\label{unicidad}
A closed Riemann surface admits at most one generalized Fermat group.
\end{theo}

In terms of Fuchsian groups, the above uniqueness result is equivalent to the following commutator rigidity property.
 
\begin{coro}\label{rigido}
For $j=1,2$, let $\Gamma_{j}<{\rm Aut}({\mathbb H}^{2})$ be a co-compact Fuchsian group with  signature $(0;k_{j},\stackrel{n_{j}+1}{\ldots},k_{j})$, where $k_{j}, n_{j} \geq 2$ and $(n_{j}-1)(k_{j}-1)>2$.
 If $\Gamma'_{1}=\Gamma'_{2}$, then $\Gamma_{1}=\Gamma_{2}$.
\end{coro}

By Theorem \ref{unicidad}, a homology Riemann surface admits at most one generalized Fermat group. One may wonder if this uniqueness property holds for general homology groups. As a consequence of the results in \cite{Singerman}, the generic homology closed Riemann surface admits only one homology group. 
In Section \ref{Sec:ejemplos}, we show explicit examples to see that uniqueness of homology groups is not always true (and, moreover, they might be either normal or non-normal subgroups). In Example \ref{ejemplo1}, the surface has genus two and it has two different conjugated homology groups isomorphic to ${\mathbb Z}_{6}$ (in particular, these homology groups are not normal subgroups). In Example \ref{ejemplo2}, the surface is hyperelliptic of genus $g \geq 2$ even, and it admits two non-isomorphic homology groups, one isomorphic to ${\mathbb Z}_{6}$ and the other isomorphic to ${\mathbb Z}_{6} \times {\mathbb Z}_{2}$, both of them being normal subgroups.

\medskip

As noted from Example \ref{ejemplo1}, a homology group $A$ of $S$ might not be a normal subgroup of ${\rm Aut}(S)$. We proceed to provide a description of the normalizers $N_{A}$ of $A$ in ${\rm Aut}(S)$. First, we need some definitions.
Let $\mu_{S,A}$ be the least common multiple of the branch orders of the conical points of $S/A$. For each $p \in S$, with a non-trivial $A$-stabilizer $A_{p}$, set $n_{p}:=\mu_{S,A}/|A_{p}|$. Let $S^{orb,A}$ be the Riemann orbifold whose underlying Riemann surface is $S$ and its cone points are those points $p \in S$ with non-trivial $A_{p}$ and $n_{p} \geq 2$ (which is the corresponding cone order). Let ${\rm Aut}(S^{orb,A})$ be the group of conformal automorphisms of $S$ keeping invariant the above cone points together their orders. One may see that  $N_{A}\leq {\rm Aut}(S^{orb,A})$ and, if all cone points of $S/A$ have the same order, then $S=S^{orb,A}$.

In general, for an abelian group (not necessarily a homology group) $A<{\rm Aut}(S)$, it might happen that $N_{A} \neq {\rm Aut}(S^{orb,A})$. For instance, if we consider Klein's surface of genus three, defined by $S:=\{[x:y:z] \in {\mathbb P}^{2}: y^{7}=xz^{4}(x-z)^{2}\}$, then  ${\rm Aut}(S) \cong {\rm PSL}_{2}(7)$ (of order $168$, the maximum possible). If $A=\langle [x:y:z] \mapsto [x:e^{2 \pi i/7}y:z] \rangle \cong {\mathbb Z}_{7}$, then $S/A$ has signature $(0;7,7,7)$. If $F_{7}: =\{[x:y:z] \in {\mathbb P}^{2}: x^{7}+y^{7}+z^{7}=0\}$ (the classical Fermat curve of degree $7$, which has genus $15$) and $B=\langle [x:y:z] \mapsto [e^{2\pi i/7}x:y:z], [x:y:z] \mapsto [x:e^{2\pi i/7}y:z]\rangle \cong {\mathbb Z}_{7}^{2}$, then $S/A$ is isomorphic as orbifold to $F_{7}/B$. So 
$A$ is not a homology group of $S$. In this case, $S^{orb,A}=S$, so ${\rm Aut}(S^{orb,A})={\rm Aut}(S) \cong {\rm PSL}_{2}(7)$. As $A$ is not a normal subgroup,  $N_{A} \neq {\rm Aut}(S^{orb,A})$. In the next, we observe that this is not the case for $A$ a homology group.

\begin{theo}\label{normalidad}
If  $A$ is a homology group of the closed Riemann surface $S$, then $N_{A}={\rm Aut}(S^{orb,A})$.
\end{theo}

%%%%%%%%%%%%%%
%%%%%%%%%%%%%%
\section{Preliminaries and known facts}\label{Sec:Prelim}

%%%%%%%%%%%%%%%%
\subsection{Riemann orbifolds}
A {\it Riemann orbifold} ${\mathcal O}$ is provided by a Riemann surface $S$, called its {\it underlying Riemann surface structure},  together a discrete collection of points, say $p_{1},p_{2},... \in S$, called its {\it cone points}, where each of these cone points $p_{j}$ has associated an integer $k_{j} \geq 2$, called its {\it cone order}. 
If $S$ is a closed Riemann surface of genus $g$ (we also say that the orbifold has genus $g$), then the number of its cone points is finite, say $p_{1},..., p_{n} \in S$ and, in this case, the tuple $(g;k_{1},...,k_{n})$ is called the {\it signature} of ${\mathcal O}$. 

A {\it conformal homeomorphism} between two Riemann orbifolds is a conformal homeomorphism between the corresponding Riemann surfaces sending cone points bijectively to cone points and preserving the cone orders. If both orbifolds are the same ${\mathcal O}$, then we talk of a  {\it conformal automorphism} of ${\mathcal O}$ and we denote by ${\rm Aut}({\mathcal O})$ its group of conformal automorphisms. 

If ${\mathcal O}$ is a Riemann orbifold and $H<{\rm Aut}({\mathcal O})$ acts discontinuously (in general it will be finite),  then the quotient ${\mathcal O}/H$ is again a Riemann orbifold. Let us denote by $\pi:{\mathcal O} \to {\mathcal O}/H$ the canonical quotient map. Let $p \in {\mathcal O}$ and $H(p)$ be its $H$-stabilizer, say of order $m \geq 1$. If $p$ is not a cone point of ${\mathcal O}$ and $m \geq 2$, then $\pi(p)$ is a cone point of ${\mathcal O}/H$ of order $m$. If $p$ is a cone point of order $n \geq 2$, then $\pi(p)$ is a cone point of ${\mathcal O}/H$ of order $mn$.

%%%%%%%%%%%%%%%
\subsection{Generalized Fermat curves}\label{Sec:GFC}
Let $k \geq 2$ and $n \geq 2$ be such that $(n-1)(k-1)>2$. Let $S$ be a generalized Fermat curve of type $(k,n)$ and let  $A \cong {\mathbb Z}_{k}^{n}$ be a generalized Fermat group of type $(k,n)$ of $S$. We may identify the quotient orbifold $S/A$ with the Riemann sphere $\widehat{\mathbb C}$ and its cone points being $\infty, 0, 1, \lambda_{1},\ldots,\lambda_{n-2}$. Let $\pi_{S,A}:S \to \widehat{\mathbb C}$ be a Galois branched covering induced by the action of $A$. Below we summarize some of the previous results on these objects.

\begin{theo}[\cite{GHL, HLKP}]
Withing the above notations, the following hold.

\begin{enumerate}
\item[(a)] $A$ is the unique generalized Fermat group of type $(k,n)$ of $S$, in particular, $A$ is a normal subgroup of ${\rm Aut}(S)$.

\item[(b)] An algebraic model for $S$ is the following non-singular projective algebraic curve (a fiber product of $(n-1)$ classical Fermat curves of degree $k$)
\begin{equation}\label{generalfermat}
C_{S}:\left\{\begin{array}{ccc}
x_{1}^{k}+x_{2}^{k}+x_{3}^{k}&=&0\\
\lambda_{1}x_{1}^{k}+x_{2}^{k}+x_{4}^{k}&=&0\\
\vdots & \vdots & \vdots\\
\lambda_{n-2}x_{1}^{k}+x_{2}^{k}+x_{n+1}^{k}&=&0
\end{array}
\right\} \subset {\mathbb P}^{n}.
\end{equation}

\item[(c)] In this algebraic model, (i) $A=\langle a_{1},\ldots,a_{n}\rangle <Aut(C_{S})$, where $a_{j}$ is multiplication of the $j \mbox{ }^{\underline th}$-coordinate by a primitive $k$-root of $1$, and (ii) the Galois branched covering map $\pi_{S,A}$, in this algebraic model, is given by $$\pi_{C_S}:C_{S} \to \widehat{\mathbb C}: [x_{1}: \cdots :x_{n+1}] \mapsto - \left(\frac{x_{2}}{x_{1}}\right)^{k}.$$

\item[(d)] If $a_{n+1}=(a_{1}\cdots a_{n})^{-1}$, then (i) every element of $A$ acting with fixed points is a power of some $a_{j}$ , $j=1,\ldots, n+1$, and every fixed point of a non-trivial power of $a_{j}$ is also a fixed point of $a_{j}$.

\end{enumerate}
\end{theo}

As a consequence of the above result, and using the fact that $S$ is uniformized by the derived subgroup of the uniformizing Fuchsian group of the orbifold $S/A$, there is a short exact sequence
$1 \to A \to {\rm Aut}(S) \stackrel{\rho}{\to} {\rm Aut}(S/A) \to 1,$
where, under our identification, ${\rm Aut}(S/A)$ is the subgroup of M\"obius transformations keeping invariant the collection 
$$\{p_{1}=\infty, p_{2}=0, p_{3}=1, p_{4}=\lambda_{1},\ldots, p_{n+1}=\lambda_{n-2}\}.$$
In the above, the surjective homomorphism $\rho$ is defined by: $\pi_{S,A} \circ a = \rho(a) \circ \pi_{S,A}$, for every $a \in A$.
If $T \in {\rm Aut}(S/A)$, then it defines a permutation $\sigma_{T} \in {\mathfrak S}_{n+1}$ of these points. If $a_{T} \in A$ is such that $\rho(a_{T})=T$, then the conjugation action of $a_{T}$ on the collection $\{a_{1},\ldots,a_{n+1}\}$ is again $\sigma_{T}$ \cite{GHL}.

\begin{rema}\label{observa1}
As a generalized Fermat pair $(S,A)$ corresponds to the derived subgroup (which is a characteristic subgroup) of a Fuchsian group $\Gamma$, such that $S/A={\mathbb H}^{2}/\Gamma$, the following lifting property holds.
Let $(S_{1},A_{1})$ and $(S_{2},A_{2})$ be two generalized Fermat pairs, both of the same type $(k,n)$. Let $\pi_{S_{j},A_{j}}:S_{j} \to S_{j}/A_{j}$ be a regular (branched) covering with deck group $A_{j}$. Then, for any biholomorphism (of orbifolds) $\psi:S_{1}/A_{1} \to S_{2}/A_{2}$ there is a biholomorphism $\eta:S_{1} \to S_{2}$ such that $\pi_{S_{2},A_{2}} \circ \eta= \psi \circ \pi_{S_{1},A_{1}}$.
\end{rema}

%%%%%%%%%%%%%%%%%%%%%%
%%%%%%%%%%%%%%%%%%%%%%
\section{Proof of Theorems \ref{unicidad} and \ref{normalidad}}\label{Sec:unicidad}
%%%%%%%%%%%%%%%%%%%%%%
\subsection{Proof of Theorem \ref{unicidad}}
Let $S$ be a homology closed Riemann surface admitting generalized Fermat groups $A \cong {\mathbb Z}_{k}^{n}$ and $B \cong {\mathbb Z}_{l}^{m}$, where $k,l\geq 2$. Let us recall that $n \geq 2$ (respectively, $m \geq 2$) is uniquely determined by the genus of $S$ and $k$ (respectively, $l$).
As there is only one generalized Fermat group of a fixed type, if $k=l$, then $A=B$. So, let us assume that $k \neq l$. 

As $A$ is a normal subgroup, the group $B$ induces an abelian group $\widetilde{B}$ of conformal automorphisms of the 
Riemann orbifold ${\mathcal O}_{A}=S/A$ of signature $(0;k,\stackrel{n+1}{\ldots},k)$ (which can be identified with the Riemann sphere). 
So $\widetilde{B}$ is either isomorphic to a cyclic group ${\mathbb Z}_{q}$, $q \geq 2$, or to the Klein group ${\mathbb Z}_{2}^{2}$. 
(i) If ${\widetilde B} \cong {\mathbb Z}_{q}$, then ${\mathcal O}_{A}/\widetilde{B}$ has signature of the form $(0;k,\stackrel{\alpha}{\ldots},k,qk,\stackrel{\beta}{\ldots},qk)$, where $\alpha \geq 1$, $\beta\in \{0,1,2\}$ and $n+1=\alpha q + \beta$. 
(ii) If ${\widetilde B} \cong {\mathbb Z}_{2}^{2}$, then ${\mathcal O}_{A}/\widetilde{B}$ has signature of the form $(0;k,\stackrel{\alpha}{\ldots},k,2k,\stackrel{\beta_{1}}{\ldots},2k,2,\stackrel{\beta_{2}}{\ldots},2)$, where $\alpha \geq 0$, $\beta_{1},\beta_{2}\in \{0,1,2,3\}$, $\beta_{1}+\beta_{2}=3$ and $n+1=4\alpha + 2\beta_{1}$.

Similarly, the group $A$ induces an abelian group $\widetilde{A}$ of conformal automorphisms of 
the Riemann orbifold ${\mathcal O}_{B}=S/B$ of signature $(0;l,\stackrel{m+1}{\ldots},l)$.
So $\widetilde{A}$ is either isomorphic to a cyclic group ${\mathbb Z}_{p}$, $p \geq 2$, or to the Klein group ${\mathbb Z}_{2}^{2}$. 
(i) If ${\widetilde A} \cong {\mathbb Z}_{p}$, then ${\mathcal O}_{B}/\widetilde{A}$ has signature of the form $(0;l,\stackrel{\widehat{\alpha}}{\ldots},l,pl,\stackrel{\widehat{\beta}}{\ldots},pl)$, where $\widehat{\alpha} \geq 1$, $\widehat{\beta}\in \{0,1,2\}$ and $m+1=\widehat{\alpha} p + \widehat{\beta}$. 
(ii) If ${\widetilde A} \cong {\mathbb Z}_{2}^{2}$, then ${\mathcal O}_{B}/\widetilde{A}$ has signature of the form $(0;l,\stackrel{\widehat{\alpha}}{\ldots},l,2l,\stackrel{\widehat{\beta_{1}}}{\ldots},2l,2\stackrel{\widehat{\beta_{2}}}{\ldots},2)$, where $\widehat{\alpha} \geq 0$, $\widehat{\beta_{1}},\widehat{\beta_{2}}\in \{0,1,2,3\}$, $\widehat{\beta_{1}}+\widehat{\beta_{2}}=3$ and $m+1=4\widehat{\alpha} + 2\widehat{\beta_{1}}$.

As ${\mathcal O}_{A}={\mathcal O}_{B}=S/\langle A,B\rangle$, the two orbifolds must have the same cone points and respective cone orders. We proceed to check this in each of the possible cases.

\medskip
\noindent
(1) If $\widetilde{A}\cong {\mathbb Z}_{p}$ and $\widetilde{B} \cong {\mathbb Z}_{q}$, then (as $k \neq l$) we must have $k=pl$ and $qk=l$, from which $pq=1$, a contradiction.

\medskip
\noindent
(2) If $\widetilde{A}\cong {\mathbb Z}_{p}$ and $\widetilde{B} \cong {\mathbb Z}_{2}^{2}$, then (as $k \neq l$) we must have that $k=pl$, $\alpha=\widehat{\beta}$, and either:
\begin{enumerate}
\item[(a)] $\beta_{2}=0$, $2k=l$, $\beta_{1}=\widehat{\alpha}$,  $n+1=4\alpha + 2\beta_{1}$ and $m+1=\widehat{\alpha}p+\widehat{\beta}$.

\item[(b)] $\beta_{1}=0$, $l=2$, $\widehat{\alpha}=\beta_{2}$, $n+1=4\alpha$ and $m+1=\widehat{\alpha}p+\widehat{\beta}$.
\end{enumerate}

In case (a), as $k=pl$ and $2k=l$, we must have $2p=1$, a contradiction.  In case (b), $k=2p$, $l=2$, $n+1=4\alpha$ and $m+1=3p+\alpha$, where $\alpha \in \{1,2\}$. The genus $g$ of $S$ has the form
$$g=1+k^{n-1}((n-1)(k-1)-2)/2, \; g=1+2^{m-1}(m-3)/2.$$

If $\alpha=1$, then it follows that $3 \times 2^{3p}=32p^{2}$, which is not possible for $p \geq 2$. If $\alpha=2$, then $2^{3p}=256 p^{6}$, which is neither possible.

\medskip
\noindent
(3) If $\widetilde{A}\cong {\mathbb Z}_{2}^{2} \cong \widetilde{B}$, then (as $k \neq l$) we have the following possibilities:
\begin{enumerate}
\item[(a)] $k=2l$, and $2k=l$, which is a contradiction as $l \geq 2$.

\item[(b)] $k=2l$, $2k=2$, a contradiction as $k \geq 2$.

\item[(c)] $k=2$, $2k=l$,$2=2l$ a contradiction, as $l \geq 2$.

\item[(d)] $k=2$, $2k=2l$, $2=l$, from which $k=l$, a contradiction.
\end{enumerate}

%%%%%%%%%%%%%%%%%%%%%%%%%%%%%%%%
\subsection{Proof of Theorem \ref{normalidad}}
Let $S$ be a closed Riemann surface of genus $g \geq 2$ and let $A<{\rm Aut}(S)$ be a homology group of $S$. Then the homology orbifold ${\mathcal O}_{S,A}=S/A$ has signature $(0;k_{1},...,k_{n+1})$, where we may assume $2 \leq k_{1} \leq k_{2} \leq \cdots \leq k_{n+1}$ (satisfying Maclachlan's condition \eqref{E1}). 

Without loss of generality, we may assume that the cone points of ${\mathcal O}_{S,A}$ are given by 
$\infty$ (of order $k_{1}$), $0$ (of order $k_{2}$), $1$ (of order $k_{3}$), $\lambda_{1}$ (of order $k_{4}$),...,$\lambda_{n-2}$ (of order $k_{n+1}$). 
Let $\mu={\rm lcm}(k_{1},...,k_{n+1})$ and 
let $P_{A}:S \to \widehat{\mathbb C}$ be a regular branched covering, with deck group $A$, whose branch values are the above cone points.

Let us consider the homology orbifold ${\mathcal O}_{S,A}^{*}$ of signature $(0;\mu,\stackrel{n+1}{\cdots},\mu)$  where the cone points are the same as for ${\mathcal O}_{S,A}$ (we have only changed the order of them). As previously observed, the homology cover of ${\mathcal O}_{S,A}^{*}$ is represented by the algebraic curve
$$C_{S,A}=\left\{\begin{array}{ccc}
x_{1}^{\mu}+x_{2}^{\mu}+x_{3}^{\mu}&=&0\\
\lambda_{1}x_{1}^{\mu}+x_{2}^{\mu}+x_{4}^{\mu}&=&0\\
\vdots & \vdots & \vdots\\
\lambda_{n-2}x_{1}^{\mu}+x_{2}^{\mu}+x_{n+1}^{\mu}&=&0
\end{array}
\right\} \subset {\mathbb P}^{n},
$$
and the corresponding homology group $H_{A} \cong {\mathbb Z}_{\mu}^{n}$ (that is, $C_{S,A}/H={\mathcal O}_{S,A}^{*}$)  is generated by the transformations $a_{j}$ (multiplication of the $j$-coordinate by $w_{\mu}=e^{2 \pi i/\mu}$). 
Let $$\pi_{A}:C_{S,A} \to \widehat{\mathbb C}: [x_{1}:\cdots:x_{n+1}] \mapsto -(x_{2}/x_{1})^{\mu},$$ which is a regular branched covering, with deck group $H_{A}$ and whose branch values are $\infty, 0, 1, \lambda_{1}, \ldots, \lambda_{n-2}$.

Let $K_{A} \lhd H_{A}$ be the subgroup generated by the elements $a_{1}^{k_{1}}$, $a_{2}^{k_{2}}$,...., $a_{n}^{k_{n}}$ and $(a_{1}\cdots a_{n})^{k_{n+1}}$. The orbifold $C_{S,A}/K_{A}$ has a Riemann surface structure $S^{*}$ admitting the Abelian group $H_{A}/K_{A}$ as group of conformal automorphisms. By the construction $H_{A}/K_{A}$ is isomorphic to $A$ and 
$S^{*}$ is a homology cover of ${\mathcal O}_{S,A}$. So, we may assume $S=S^{*}$ and $A=H_{A}/K_{A}$. 

Let $Q_{A}:C_{S,A} \to S$ be a regular branched covering, with deck group $K_{A}$, such that $\pi_{A}=P_{A} \circ Q_{A}$. In this case, $S^{orb,A}=C_{S,A}/K_{A}$.
Note that the subgroup $K_{A}$ is uniquely determined by the branch values of the cone points of the orbifold ${\mathcal O}_{S,A}$.

Let $\Gamma$ be a Fuchsian group such that ${\mathbb H}^{2}/\Gamma={\mathcal O}_{S,A}^{*}$, so $C_{S,A}={\mathbb H}^{2}/\Gamma'$. By the uniqueness of $K_{A}$, there is a unique subgroup $K$ of $\Gamma$, containing $\Gamma'$ such that ${\mathbb H}^{2}/K=C_{S,A}/K_{A}$. We observe that $\Gamma'$ is the smallest normal subgroup $U$ of $K$ such that $K/U$ is isomorphic to $K_A$.
Let $\phi \in {\rm Aut}(C_{S,A}/K_A)={\rm Aut}(S^{orb,A})$ and let $\eta \in {\rm Aut}({\mathbb H}^{2})$ be a lifting of $\phi$ (so it normalizes $K$). By the uniqueness of $\Gamma'$ in $K$, $\eta$ keeps invariant it, so it descends to an automorphisms $\psi$ of $C_{S,A}$. By the uniqueness of the generalized Fermat group $H_{A}=\Gamma/\Gamma'$ \cite{HLKP}, it is also invariant under conjugation by $\psi$. As $\psi$ is a lifting of $\phi$,  $K_A$ is also invariant under $\psi$. It follows that $\phi$ normalizes $A$.

%%%%%%%%%%%%%%%%%%%%
%%%%%%%%%%%%%%%%%%%%
\section{Examples}\label{Sec:ejemplos}

%%%%%%%%%%%%%%%%%%%
\begin{example}
Let $S$ be the genus three hyperelliptic Riemann surface defined by the hyperelliptic curve $w^{2}=u^{8}-1.$
This surface admits the following automorphisms:
$$a(u,w)=\left( \sqrt{i} \; u,w \right),\; b(u,w)=(-u,-w), \; c(u,w)=\left(1/u, iw/u^{4}\right),$$
of respective orders $8$, $2$ and $4$. If $A=\langle a,b\rangle \cong {\mathbb Z}_{8} \times {\mathbb Z}_{2}$, then $S/A$ has signature $(0;2,8,8)$. It can be checked that $A$ is a homology group of $S$. In this case, the points in $S$ projecting to those of order $2$ in $S/A$ are the eight Weierstrass points and these are the cone points of the orbifold $S^{orb,A}$, each of them with cone order $4$, so ${\rm Aut}(S^{orb,A})={\rm Aut}(S)$. Theorem \ref{normalidad} asserts  that $A$ is a normal subgroup of ${\rm Aut}(S)$.
\end{example}

%%%%%%%%%%%%%%%%%%%%
\begin{example}\label{ejemplo1}
Let $S$ be the genus two Riemann surface defined by the hyperelliptic curve $y^{2}=x(x^{4}-1)$. This surface admits the order six automorphism
$$\eta(x,y)=\left( i (1+x)/(1-x), 2(1-i)y/(x-1)^{3}\right)$$
and $\eta^{3}(x,y)=(x,-y)$ is the hyperelliptic involution. If $A=\langle \eta \rangle \cong {\mathbb Z}_{6}$, then one may see that the quotient orbifold $S/A$ has signature $(0;2,2,3,3)$ and $A$ is a homology group of $S$. On $S$ we also have the order four automorphism $\rho(x,y)=(-x,iy)$. Then $B=\rho A \rho^{-1} \cong {\mathbb Z}_{6}$ is also a homology group of $S$. As
$$\rho \circ \eta \circ \rho^{-1}(x,y)=\left( -i(1-x)/(1+x),-2(1-i)y/(1+x)^{3}\right) \notin A,$$
we see that $A \neq B$. In particular, $S$ has two different homology groups, both isomorphic to ${\mathbb Z}_{6}$, and $N_{A} \neq {\rm Aut}(S)$.
\end{example}

%%%%%%%%%%%%%%%%%%%%
\begin{example}\label{ejemplo2}
Let $g \geq 2$ be an even integer and let $S$ be the genus $g$ Riemann surface  defined by the hyperelliptic curve $y^{2}=x^{2g+2}-1$. This Riemann surface admits the order $2g+2$ automorphism
$\alpha(x,y)=(e^{\pi i/(g+1)}x,y)$
and $\tau(x,y)=(x,-y)$ is its hyperelliptic involution. If $A=\langle \alpha, \tau \rangle \cong {\mathbb Z}_{2g+2} \times {\mathbb Z}_{2}$, then one may see that the quotient orbifold $S/A$ has signature $(0;2,2g+2,2g+2)$ and $A$ is a homology group of $S$. Similarly, if 
$B=\langle \alpha^{2}, \tau \rangle \cong {\mathbb Z}_{2g+2}$, then $S/B$ has signature $(0;2,2,g+1,g+1)$ and $B$ is again a homology group of $S$.
In particular, $S$ has two different homology groups, one isomorphic to ${\mathbb Z}_{2g+2}$ and the other to ${\mathbb Z}_{2g+2} \times {\mathbb Z}_{2}$. Note that, by Theorem \ref{normalidad},  $N_{A}={\rm Aut}(S)$ as the orbifold points of $S^{orb,A}$ are exactly the Weierstrass points, each one with cone order $g+1$; in particular, $A$ is a normal subgroup. Similarly, it can be seen that $B$ is also a normal subgroup.
\end{example}

%%%%%%%%%%%%%%
%%%%%%%%%%%%%%

\end{document}